# Nonresonance and global existence of prestressed nonlinear elastic waves

By Thomas C. Sideris*

## 1. Introduction

This article considers the existence of global classical solutions to the Cauchy problem in nonlinear elastodynamics. The unbounded elastic medium is assumed to be homogeneous, isotropic, and hyperelastic. As in the theory of 3D nonlinear wave equations in three space dimensions, global existence hinges on two basic assumptions. First, the initial deformation must be a *small displacement* from equilibrium, in this case a prestressed homogeneous dilation of the reference configuration, and equally important, the nonlinear terms must obey a type of *nonresonance* or null condition. The omission of either of these assumptions can lead to the breakdown of solutions in finite time. In particular, nonresonance complements the genuine nonlinearity condition of F. John, under which arbitrarily small spherically symmetric displacements develop singularities (although one expects this to carry over to the nonsymmetric case, as well), [4]. John also showed that small solutions exist almost globally [5] (see also [10]). Formation of singularities for large displacements was illustrated by Tahvildar-Zadeh [16].

The nonresonance condition introduced here represents a substantial improvement over our previous work on this topic [13]. To explain the difference roughly, our earlier version of the null condition forced the cancellation of all nonlinear wave interactions to first order along the characteristic cones. Here, only the cancellation of nonlinear wave interactions among individual wave families is required. The difficulty in realizing this weaker version is that the decomposition of elastic waves into their longitudinal and transverse components involves the nonlocal Helmholtz projection, which is ill-suited to nonlinear analysis. However, our decay estimates make clear that only the leading contribution of the resonant interactions *along the characteristic cones* is potentially dangerous, and this permits the usage of approximate *local* decompositions.

*Supported in part by the National Science Foundation.

The limitations of the earlier version of the null condition were exposed in the work of Tahvildar-Zadeh who first considered small perturbations from an arbitrary prestressed homogeneous dilation [16]. Using an expansion based on small shear strains, he observed that the null condition in [13] placed physically unrealistic restrictions on the growth of the stored energy function for extreme values of the dilational strain. Our revised null condition avoids this defect. In fact, an expansion involving small distortional strain (which to first order is equivalent with small shear strain) reveals that the dominant behavior of the stored energy function for arbitrary dilational strain is determined by the propagation speeds of the medium. The null condition enters as a higher order correction.

Outside of the inevitable energy estimates, the techniques of the existence proof represent an evolution of the ideas initiated in [10] and expanded in [13]. The central hurdle to be overcome is the fact that the equations of motion for elasticity are not Lorentz invariant. For the nonlinear wave equation, a well-oiled machine based on Lorentz invariance, the so-called generalized energy method, has been developed by Klainerman for the construction of solutions [8]. Its attractive feature is the avoidance of direct estimation of the fundamental solution of the wave equation through the use of weighted Sobolev inequalities involving the generators of the Lorentz group. The generators are effective because they commute with the wave operator. In the case of elasticity, it is still possible to get weighted estimates with respect to a smaller number of generators, but the resulting estimates are weaker. In particular, the best decay is available only for second (or higher) derivatives of the solution. In our earlier work [10], [13], this was dealt with by smoothing and the consequent introduction of a somewhat unnatural and cumbersome nonlocal energy. The proof presented here represents a refinement also in that the use of the nonlocal energy is now completely avoided. Given the significant differences from our earlier work, we have endeavored to make a self-contained presentation.

The results presented herein were announced in [14]. R. Agemi has recently completed a manuscript dealing with the existence of solutions near the unstressed reference configuration under the same null condition used here [1]. His proof relies on direct estimation of the fundamental solution. The approach taken here can be adapted to the case of a system of coupled quadratically nonlinear wave equations in 3D with multiple propagation speeds [15].

The remainder of this introduction will be devoted to the description of the basic notation and the formulation of the equations of motion, leading as speedily as possible to a statement of the global existence theorem. In the following section, we then explore the origins of our assumptions in an effort to convince the reader of their transparency. The rest of the paper is devoted to the proof of global existence. The key components: commutation, null form estimates, weighted $L^\infty$ estimates, and weighted $L^2$ estimates, are assembled along the way in a series of four propositions.

1.1. *Notation.* Partial derivatives will be written as

$$\partial_0 = \partial_t = \frac{\partial}{\partial t} \quad \text{and} \quad \partial_i = \frac{\partial}{\partial x^i}, \quad i = 1, 2, 3.$$

We will also abbreviate

$$\partial = (\partial_0, \partial_1, \partial_2, \partial_3), \quad \text{and} \quad \nabla = (\partial_1, \partial_2, \partial_3).$$

The so-called angular momentum operators are the vector fields

$$\Omega = (\Omega_1, \Omega_2, \Omega_3) = x \wedge \nabla,$$

$\wedge$ being the usual vector cross product. Crucially, the spatial partial derivatives can be decomposed into radial and angular components

(1.1) $$\nabla = \frac{x}{r}\partial_r - \frac{x}{r^2} \wedge \Omega, \quad \text{where} \quad r = |x|, \quad \text{and} \quad \partial_r = \frac{x}{r} \cdot \nabla.$$

A central role is played by the generators of simultaneous rotations which will be seen in Section 3.1 to arise naturally in connection with the symmetries of isotropic materials. They are given by the following vectorial modification of the angular momentum operators

(1.2a) $$\tilde{\Omega}_\ell = \Omega_\ell I + U_\ell,$$

with

(1.2b) $$U_1 = \begin{bmatrix} 0 & 0 & 0 \\ 0 & 0 & 1 \\ 0 & -1 & 0 \end{bmatrix}, \quad U_2 = \begin{bmatrix} 0 & 0 & -1 \\ 0 & 0 & 0 \\ 1 & 0 & 0 \end{bmatrix}, \quad U_3 = \begin{bmatrix} 0 & 1 & 0 \\ -1 & 0 & 0 \\ 0 & 0 & 0 \end{bmatrix}.$$

Equally important will be the scaling operator

$$S = t\partial_t + r\partial_r,$$

but we also shall show in Section 3.1 why it is more precise to use the generator of changes of scale

(1.3) $$\tilde{S} = S - 1.$$

The eight vector fields will be written as $\Gamma = (\Gamma_0, \ldots, \Gamma_7) = (\partial, \tilde{\Omega}, \tilde{S})$. The commutator of any two $\Gamma$'s is either 0 or is in the span of $\Gamma$, and in particular, the commutator of $\nabla$ and $\Gamma$ is in the span of $\nabla$. Schematically, we write

(1.4) $$[\nabla, \Gamma] = \nabla.$$

By $\Gamma^a$, $a = (a_1, \ldots, a_\kappa)$, we denote an ordered product of $\kappa = |a|$ vector fields $\Gamma_{a_1} \cdots \Gamma_{a_\kappa}$.[1]

---

[1] Note that this differs from the standard multi-index notation.

The energy associated to the linearized operator to be defined in Section 1.2 is

$$E_1(u(t)) = \tfrac{1}{2} \int_{\mathbb{R}^3} \left[ |\partial_t u(t)|^2 + c_2^2 |\nabla u(t)|^2 + (c_1^2 - c_2^2)(\nabla \cdot u(t))^2 \right] dx,$$

and higher order energies are defined through

$$E_\kappa(u(t)) = \sum_{|a| \leq \kappa-1} E_1(\Gamma^a u(t)).$$

In order to describe the solution space we also introduce the time-independent analog of $\Gamma$. Set

$$\Lambda = (\Lambda_1, \ldots, \Lambda_7) = (\nabla, \widetilde{\Omega}, r\partial_r - 1).$$

Then the $\Lambda$'s have the same commutation properties as the $\Gamma$'s. Define

$$H_\Lambda^\kappa = \{ f \in L^2(\mathbb{R}^3)^3 : \Lambda^a f \in L^2(\mathbb{R}^3)^3, \ |a| \leq \kappa \},$$

with the norm

$$\|f\|_{H_\Lambda^\kappa}^2 = \sum_{|a| \leq \kappa} \|\Lambda^a f\|_{L^2(\mathbb{R}^3)}^2.$$

The solution will be constructed in the space $\dot{H}_\Gamma^\kappa(T)$ obtained by closing the set $C^\infty([0,T]; C_0^\infty(\mathbb{R}^3)^3)$ in the norm $\sup_{0 \leq t < T} E_\kappa^{1/2}(u(t))$. Thus,

$$\dot{H}_\Gamma^\kappa(T) \subset \left\{ u(t,x) : \partial u(t, \cdot) \in \bigcap_{j=0}^{\kappa-1} C^j([0,T]; H_\Lambda^{\kappa-1-j}) \right\}.$$

By (3.20a), it will follow that $\dot{H}_\Gamma^\kappa(T) \subset C^{\kappa-2}([0,T] \times \mathbb{R}^3)$.

We define the orthogonal projections onto radial and transverse directions by

(1.5a) $$P_1 u(x) = \frac{x}{r} \otimes \frac{x}{r} u(x) = \frac{x}{r} \left\langle \frac{x}{r}, u(x) \right\rangle$$

and

(1.5b) $$P_2 u(x) = [I - P_1] u(x) = -\frac{x}{r} \wedge \left( \frac{x}{r} \wedge u(x) \right).$$

Employing the notation $\langle \rho \rangle = (1 + |\rho|^2)^{1/2}$, we will use the following weighted $L^2$-norm

(1.6) $$\mathcal{X}_\kappa(u(t)) = \sum_{\alpha=1}^{2} \sum_{\beta=0}^{3} \sum_{\ell=1}^{3} \sum_{|a| \leq \kappa-2} \|\langle c_\alpha t - r \rangle P_\alpha \partial_\beta \partial_\ell \Gamma^a u(t)\|_{L^2}.$$

1.2. *The equations of motion.* Consider a homogeneous elastic material filling space. Assume that its density in its undeformed state is unity. The fundamental unknown is the deformation $\varphi : \mathbb{R} \times \mathbb{R}^3 \to \mathbb{R}^3$ which is an orientation-preserving diffeomorphism taking a material point $x \in \mathbb{R}^3$ in the reference configuration to its position $\varphi(t, x) \in \mathbb{R}^3$ at time $t$. Consequently, the deformation gradient $F = \nabla \varphi$ must satisfy $\det F > 0$.

A hyperelastic material is one for which there exists a stored energy function $\sigma(F)$, representing the potential energy term in the Hamiltonian. The Piola-Kirchhoff stress tensor is given by $\partial \sigma(F)/\partial F$. A material is frame indifferent, respectively, isotropic if the conditions

(1.7) $$\sigma(QF) = \sigma(F) \quad \text{and} \quad \sigma(FQ) = \sigma(F)$$

hold for every proper orthogonal matrix $Q$. It is well-known that (1.7) implies that the stored energy function $\sigma$ depends on $F$ only through the principal invariants of the strain tensor $F^T F$; see [11, p. 192].

Under these assumptions, the equations of motion are[2]

(1.8) $$\frac{\partial^2 \varphi^i}{\partial t^2} - \frac{\partial}{\partial x^\ell}\left(\frac{\partial \sigma}{\partial F_\ell^i}(\nabla \varphi)\right) = 0.$$

We shall consider deformations of the form $\varphi(t,x) = \lambda x + u(t,x)$, $\lambda > 0$, in which $u(t,x)$ represents a small displacement from a homogeneous dilation, so that, in particular, $\varphi(t,\cdot)$ is a bijection. For small displacements $u$, the standard linear theory outlined in Section 2.2 ensures that the system (1.8) is hyperbolic. The long-time behavior of solutions of quasilinear wave equations in 3D is determined by the structure of the quadratic portion on the nonlinearity [2], [3], [7], [9], [12], and so there is no essential loss of generality in considering the truncated version of (1.8) obtained by expansion to second order about the equilibrium. This yields the quasilinear system

(1.9a) $$Lu \equiv \partial_t^2 u - Au = N(u,u),$$

with

(1.9b) $$\begin{cases} (Au)^i = A_{\ell m}^{ij}(\lambda) \partial_\ell \partial_m u^j, \\ A_{\ell m}^{ij}(\lambda) = \dfrac{\partial^2 \sigma}{\partial F_\ell^i \partial F_m^j}(\lambda I), \end{cases}$$

and

(1.9c) $$\begin{cases} N(u,v)^i = B_{\ell m n}^{ijk}(\lambda) \partial_\ell\left(\partial_m u^j \partial_n v^k\right), \\ B_{\ell m n}^{ijk}(\lambda) = \dfrac{\partial^3 \sigma}{\partial F_\ell^i \partial F_m^j \partial F_n^k}(\lambda I). \end{cases}$$

When it is secondary, the dependence of the coefficients on the parameter $\lambda$ will be suppressed. It follows from the definitions (1.9b), (1.9c) that the coefficients are symmetric with respect to pairs of indices

(1.10) $$A_{\ell m}^{ij} = A_{m\ell}^{ji} \quad \text{and} \quad B_{\ell m n}^{ijk} = B_{m\ell n}^{jik} = B_{\ell n m}^{ikj}.$$

---

[2] Repeated indices are always summed, regardless of their position up or down.

Thus, the form (1.9c) is symmetric and

(1.11) $$N(u,v) = N(v,u).$$

1.3. *The global existence theorem.* The basic assumptions of linear elasticity, summarized in Section 2.2, ensure the hyperbolicity of the linear operator. It will be seen in (2.6b) to have the form

(H) $\quad L = (\partial_t^2 - c_2^2(\lambda)\Delta)I - (c_1^2(\lambda) - c_2^2(\lambda))\nabla \otimes \nabla, \quad c_1^2(\lambda) > \frac{4}{3}c_2^2(\lambda) > 0.$

As will be more fully explained in Section 2.3, the null condition in this instance is:[3]

(N) $\quad B^{ijk}_{\ell mn}(\lambda)\xi_i\xi_j\xi_k\xi_\ell\xi_m\xi_n = 0, \quad \text{for all } \xi \in S^2.$

With these essentials, we can now give a precise statement of the result to be shown in Section 3.

THEOREM 1.1. *For every value of $\lambda \in \mathbb{R}^+$ such that the hyperbolicity condition* (H) *and the null condition* (N) *hold, the initial value problem for* (1.9a)–(1.9c) *with initial data*

$$\partial u(0) \in H^{\kappa-1}_\Lambda, \quad \kappa \geq 9,$$

*has a unique global solution* $u \in \dot{H}^\kappa_\Gamma(T)$ *for every $T > 0$, provided that*

(1.12) $\quad E_{\kappa-2}(u(0)) \, \exp\left[C(\lambda)E_\kappa^{1/2}(u(0))\right] \leq \varepsilon(\lambda),$

*and $\varepsilon(\lambda)$ is sufficiently small, depending on $\lambda$. The solution satisfies the bounds*

$$E_{\kappa-2}(u(t)) \leq \varepsilon(\lambda) \quad \text{and} \quad E_\kappa(u(t)) \leq 2E_\kappa(u(0))\langle t \rangle^{C(\lambda)\sqrt{\varepsilon(\lambda)}},$$

*for all $t \geq 0$.*

## 2. Preliminaries

2.1. *A word about invariants.* In addition to their definitions as the elementary symmetric functions of the eigenvalues, the invariants $\mathcal{I}(C)$ of $3 \times 3$ matrix $C$ are conveniently expressed as

(2.1a) $\quad \mathcal{I}_1(C) = \text{tr } C, \quad \mathcal{I}_2(C) = \frac{1}{2}[(\text{tr } C)^2 - \text{tr } C^2], \quad \mathcal{I}_3(C) = \det C.$

The invariants of any two matrices which differ by a multiple of the identity are linearly related. Thus, if $C' = \mathbf{z}\, I + C$, then

(2.1b) $\begin{aligned} \mathcal{I}_1(C') &= 3\mathbf{z} + \mathcal{I}_1(C) \\ \mathcal{I}_2(C') &= 3\mathbf{z}^2 + 2\mathbf{z}\,\mathcal{I}_1(C) + \mathcal{I}_2(C) \\ \mathcal{I}_3(C') &= \mathbf{z}^3 + \mathbf{z}^2\,\mathcal{I}_1(C) + \mathbf{z}\,\mathcal{I}_2(C) + \mathcal{I}_3(C), \end{aligned}$

as can easily be verified by comparing eigenvalues.

---

[3]Our earlier version was $B^{ijk}_{\ell mn}(\lambda)\xi_\ell\xi_m\xi_n = 0$, for all $\xi \in S^2$ and all $i,j,k$.

In order to make sense of the coming assumptions about the stored energy function, it will be necessary to explore its dependence on the invariants of the strain matrix. When considering perturbations from an equilibrium, it is natural to introduce the invariants of the perturbation. We shall alternate between the strain matrix with which computations are most easily performed and its square root, the stretch matrix, in terms of which the results are most concisely expressed. The following table is meant to collect the notation to be used:

$$
\begin{array}{ll}
\mathcal{I}_k(F^T F) = \mathtt{i}_k & \mathcal{I}_k(F^T F - \lambda^2 I) = \mathtt{j}_k \\
\mathcal{I}_k(\sqrt{F^T F}) = \mathtt{r}_k & \mathcal{I}_k(\sqrt{F^T F} - \lambda I) = \mathtt{s}_k
\end{array}
\qquad k = 1, 2, 3.
\tag{2.2a}
$$

Thus, by (2.1b), $\mathtt{i}$ and $\mathtt{r}$ are linearly related to $\mathtt{j}$ and $\mathtt{s}$, respectively.

Through comparison of eigenvalues, the invariants $\mathtt{i}$ and $\mathtt{r}$ are seen to satisfy

$$\mathtt{i}_1 = \mathtt{r}_1^2 - 2\mathtt{r}_2, \quad \mathtt{i}_2 = \mathtt{r}_2^2 - 2\mathtt{r}_1\mathtt{r}_3, \quad \mathtt{i}_3 = \mathtt{r}_3^2. \tag{2.2b}$$

This transformation is invertible in a neighborhood of the equilibrium, by the implicit function theorem. From (2.2b), a short computation using (2.1b) produces

$$
\begin{aligned}
\mathtt{j}_1 &= 2\lambda \mathtt{s}_1 + \mathtt{s}_1^2 - 2\mathtt{s}_2, \\
\mathtt{j}_2 &= 4\lambda^2 \mathtt{s}_2 + 2\lambda \mathtt{s}_1 \mathtt{s}_2 - 6\lambda \mathtt{s}_3 + \mathtt{s}_2^2 - 2\mathtt{s}_1 \mathtt{s}_3, \\
\mathtt{j}_3 &= 8\lambda^3 \mathtt{s}_3 + 4\lambda^2 \mathtt{s}_1 \mathtt{s}_3 + 2\lambda \mathtt{s}_2 \mathtt{s}_3 + \mathtt{s}_3^2,
\end{aligned}
\tag{2.2c}
$$

which is also invertible near the origin.

*Note.* As a consequence of these formulae, it is possible to switch between any of the sets of variables $\mathtt{i}$, $\mathtt{j}$, $\mathtt{r}$, $\mathtt{s}$. In particular, by (2.1b), (2.2c), we can write

$$\sigma(F) = \sigma(\mathtt{i}) = \tau(\lambda, \mathtt{s}). \tag{2.2d}$$

**2.2. The linear operator.** In order that the linear operator $L$ defined in (1.9a) be hyperbolic, we need to impose an ellipticity condition for the operator $A$ in (1.9b).

LEMMA 2.1. *Let $A$ be the operator defined in (1.9b) from an isotropic stored energy function $\sigma$. The following statements are equivalent:*
- *$A$ is elliptic.*
- *The symbol of $A$, $A^{ij}(\xi) \equiv A^{ij}_{\ell m}(\lambda)\xi_\ell \xi_m$, is positive definite:*

$$A^{ij}_{\ell m}(\lambda)\xi_\ell \xi_m x^i x^j > 0, \quad \text{for all } x, \xi \in S^2. \tag{2.3a}$$

- *The Legendre-Hadamard condition is satisfied:*

$$D_\varepsilon^2 \sigma(\lambda I + \varepsilon\, x \otimes \xi)\Big|_{\varepsilon=0} > 0, \quad \text{for all } x, \xi \in S^2. \tag{2.3b}$$

- *The positivity conditions*

(2.3c) $$\tau_{11}(\lambda, 0) > 0 \quad \text{and} \quad \tau_1(\lambda, 0) - \lambda \tau_2(\lambda, 0) > 0$$

*hold for the stored energy function $\tau$ in* (2.2d).[4]

*Proof.* Condition (2.3a) is merely the definition of ellipticity.

Thanks to (1.9b), the derivative in (2.3b) is equal to the expression in (2.3a), so the two conditions are equivalent.

From (2.1a), it is a simple matter to check that when $F = \lambda I + \varepsilon\, x \otimes \xi$, $x, \xi \in S^2$, the invariants of $F^T F - \lambda^2 I$ satisfy

(2.4a) $$\mathsf{j}_1 = 2\lambda\varepsilon\langle x, \xi\rangle + \varepsilon^2, \quad \mathsf{j}_2 = -\lambda^2\varepsilon^2 |x \wedge \xi|^2, \quad \mathsf{j}_3 = 0.$$

Combining (2.4a) and (2.2c), we find that $\mathsf{s}_3 = 0$ and

(2.4b) $$D_\varepsilon \mathsf{s}_1|_{\varepsilon=0} = \langle x, \xi\rangle \qquad D_\varepsilon \mathsf{s}_2|_{\varepsilon=0} = 0$$

$$D_\varepsilon^2 \mathsf{s}_1|_{\varepsilon=0} = \tfrac{1}{2\lambda}|x \wedge \xi|^2 \qquad D_\varepsilon^2 \mathsf{s}_2|_{\varepsilon=0} = -\tfrac{1}{2}|x \wedge \xi|^2.$$

So with (2.2d), (2.4b), we can explicitly compute the derivative in (2.3b):

(2.4c) $$\begin{aligned}D_\varepsilon^2 \sigma(\lambda I + \varepsilon\, x \otimes \xi)|_{\varepsilon=0} &= D_\varepsilon^2 \tau(\lambda, \mathsf{s}_1, \mathsf{s}_2, 0)|_{\varepsilon=0} \\ &= \tau_{11}(\lambda, 0)\, \langle x, \xi\rangle^2 + \tfrac{1}{2\lambda}[\tau_1(\lambda, 0) - \lambda\tau_2(\lambda, 0)]\,|x \wedge \xi|^2,\end{aligned}$$

thereby showing the equivalence of (2.3c) and (2.3b). □

We shall now impose conditions (2.3a)–(2.3c). Considering (2.3c), we define positive constants $c_1(\lambda)$ and $c_2(\lambda)$, representing the propagation speeds, by

(2.5a) $$c_1^2(\lambda) = \tau_{11}(\lambda, 0), \qquad c_2^2(\lambda) = \tfrac{1}{2\lambda}[\tau_1(\lambda, 0) - \lambda \tau_2(\lambda, 0)].$$

Following the standard linear theory, we shall assume that the reference configuration, $\lambda = 1$, $\mathsf{s} = 0$, is a stress-free state

(2.5b) $$\tau_1(1, 0) = 0.$$

It is necessary that the speeds remain distinct, and to be consistent with linear theory we assume that the bulk and shear moduli are positive:

(2.5c) $$c_1^2(\lambda) - \tfrac{4}{3}c_2^2(\lambda) > 0 \quad \text{and} \quad c_2^2(\lambda) > 0.$$

From (1.10), (2.3a), (2.4c), and (2.5a) we see that the symbol matrix $A(\xi)$ is positive definite and symmetric with the squared speeds as its eigenvalues. The corresponding eigenspaces are the one-dimensional span of the unit

---

[4] Here and later on, subscripts for $\tau$ indicate derivatives in $\mathsf{s}$.

vector $\xi$ and its orthogonal complement. Consequently, we can write the symbol in spectral form

$$(2.6a) \qquad A(\xi) = c_1^2(\lambda)\xi \otimes \xi + c_2^2(\lambda)[I - \xi \otimes \xi], \quad \xi \in S^2,$$

or in other words,

$$(2.6b) \qquad A = c_2^2(\lambda)\Delta I + (c_1^2(\lambda) - c_2^2(\lambda))\nabla \otimes \nabla.$$

For each direction $\xi \in S^2$, we have two families of elementary plane wave solutions of $Lu = 0$, namely

$$(2.7) \quad \begin{aligned} \mathcal{W}_1(\xi) &= \{\alpha\xi \exp \mathbf{i}\beta[\langle x, \xi\rangle - c_1(\lambda)t] : \alpha, \beta \in \mathbb{R}\}, \\ \mathcal{W}_2(\xi) &= \{\eta \exp \mathbf{i}\beta[\langle x, \xi\rangle - c_2(\lambda)t] : \langle \eta, \xi\rangle = 0, \beta \in \mathbb{R}\}. \end{aligned}$$

These elementary solutions represent longitudinal and transverse waves propagating in the direction $\xi$ with speeds $c_1(\lambda)$ and $c_2(\lambda)$, respectively.

2.3. *The nonlinearity.* Global existence requires a further nonresonance condition, the so-called null condition, linking the quadratic portion of the nonlinearity with the linear operator: the quadratic interaction of elementary waves of each wave family only produces waves in the other family.[5] This idea is expressed in the following:

*Definition* 2.1. The quadratic nonlinearity $N$ defined in (1.9c) is null with respect to the linear operator $L$ defined in (1.9a) and (1.9b) provided that

$$(2.8a) \qquad \langle u, N(v, w)\rangle = 0,$$

for all resonant triples

$$(u, v, w) \in \mathcal{W}_\alpha(\xi) \times \mathcal{W}_\alpha(\xi) \times \mathcal{W}_\alpha(\xi), \quad \alpha = 1, 2.$$

In terms of the coefficients $B^{ijk}_{\ell mn}$ of the nonlinearity, we see by direct substitution that condition (2.8a) is equivalent to

$$(2.8b) \qquad B^{ijk}_{\ell mn}(\lambda)\xi_i\xi_j\xi_k\xi_\ell\xi_m\xi_n = 0, \quad \text{for all } \xi \in S^2,$$

$$(2.8c) \qquad B^{ijk}_{\ell mn}(\lambda)\eta^{(1)}_i\eta^{(2)}_j\eta^{(3)}_k\xi_\ell\xi_m\xi_n = 0, \quad \begin{aligned}&\text{for all } \xi, \eta^{(a)} \in S^2 \\ &\text{with } \langle\xi, \eta^{(a)}\rangle = 0.\end{aligned}$$

The following result relates this condition to the stored energy function and the stress tensor. In particular, it shows that condition (2.8c) for the transverse waves is redundant in the isotropic case.

LEMMA 2.2. *Assume that the nonlinear quadratic form $N$ defined in (1.9c) arises from an isotropic stored energy function $\sigma$. The following statements are equivalent*:

---

[5]A connection between plane waves and the null condition was first noted in [6].

- *The nonlinear form N is null with respect to L.*
- *Condition (2.8b) holds.*
- *The stored energy function satisfies*

$$\text{(2.9a)} \qquad D_\varepsilon^3 \sigma(\lambda I + \varepsilon\, \xi \otimes \xi)|_{\varepsilon=0} = 0, \quad \text{for all } \xi \in S^2.$$

- *The Piola-Kirchhoff stress tensor $\Sigma(F) = \partial\sigma(F)/\partial F$ satisfies*

$$\text{(2.9b)} \qquad \operatorname{tr} D_\varepsilon^2 \Sigma(\lambda I + \varepsilon\, \xi \otimes \xi)|_{\varepsilon=0}\, \xi \otimes \xi = 0, \quad \text{for all } \xi \in S^2.$$

- *The degeneracy condition*

$$\text{(2.9c)} \qquad \tau_{111}(\lambda, 0) = 0$$

*holds for $\tau$ in (2.2d).*

*Proof.* It is clear that conditions (2.8b), (2.9a), and (2.9b) are equivalent, while condition (2.8c) is equivalent to

$$\text{(2.10a)} \quad D_{\varepsilon_1} D_{\varepsilon_2} D_{\varepsilon_3} \sigma(\lambda I + (\varepsilon_1 \eta^{(1)} + \varepsilon_2 \eta^{(2)} + \varepsilon_3 \eta^{(3)}) \otimes \xi)|_{\varepsilon_a=0} = 0,$$
$$\text{for all } \xi, \eta^{(a)} \in S^2 \text{ with } \langle \xi, \eta^{(a)} \rangle = 0.$$

To see that this holds for all isotropic materials, take

$$F = \lambda I + (\varepsilon_1 \eta^{(1)} + \varepsilon_2 \eta^{(2)} + \varepsilon_3 \eta^{(3)}) \otimes \xi$$

with $\xi, \eta^{(a)} \in S^2$ and $\langle \xi, \eta^{(a)} \rangle = 0$. It follows from (2.4a) with $\varepsilon = |\varepsilon_a \eta^{(a)}|$ and $x = \varepsilon_a \eta^{(a)}/\varepsilon$ that now

$$\text{(2.10b)} \qquad \begin{aligned}
\mathtt{j}_1 &= |\varepsilon_1 \eta^{(1)} + \varepsilon_2 \eta^{(2)} + \varepsilon_3 \eta^{(3)}|^2, \\
\mathtt{j}_2 &= -\lambda^2 |(\varepsilon_1 \eta^{(1)} + \varepsilon_2 \eta^{(2)} + \varepsilon_3 \eta^{(3)}) \wedge \xi|^2, \\
\mathtt{j}_3 &= 0.
\end{aligned}$$

Note that (2.10b) is quadratic in $\varepsilon_a$. Therefore, regarding $\sigma(\mathtt{i})$ as a function of $\mathtt{j}$ through (2.1b), we see that (2.10a) is true without further assumptions.

Finally, to get (2.9c) let us consider (2.9a). Setting $F = \lambda I + \varepsilon\, \xi \otimes \xi$, we get from (2.4a) that $\mathtt{j}_1 = 2\lambda\varepsilon + \varepsilon^2$, $\mathtt{j}_2 = \mathtt{j}_3 = 0$, and so from (2.2c) it follows that $\mathtt{s}_1 = \varepsilon$, $\mathtt{s}_2 = \mathtt{s}_3 = 0$. Switching variables, as in (2.2d), we see that $\sigma(\lambda I + \varepsilon\, \xi \otimes \xi) = \tau(\lambda, \varepsilon, 0, 0)$, and the null condition (2.9a) for the longitudinal waves reduces to (2.9c). □

*Note.* When $\lambda = 1$, the condition (2.9c) is complementary to John's *genuine nonlinearity* condition $\tau_{111}(1, 0) > 0$ which leads to formation of singularities in small spherically symmetric displacements [4].

*Note.* In terms of the Cauchy stress tensor, $\mathtt{T}(F) = \det(F)^{-1} \Sigma(F) F^T$, the null condition can be characterized as

$$\operatorname{tr} D_\varepsilon^2 \mathtt{T}(\lambda I + \varepsilon\, \xi \otimes \xi)|_{\varepsilon=0}\, \xi \otimes \xi = 0, \quad \text{for all } \xi \in S^2.$$

2.4. *Deciphering the conditions.* This section will illustrate that our conditions can be satisfied for all values of $\lambda$ with physically realistic choices of the stored energy function.

An expansion of the stored energy function can be based on the distortional strain matrix $C = \sqrt{F^T F} - \frac{1}{3}\mathbf{r}_1 I$, where, with the notation in (2.2a), $\mathbf{r}_1 = \operatorname{tr} \sqrt{F^T F}$ is the dilational strain. (The distortional strain is essentially the linearization of the shear strain considered in [16].) Notice that $\operatorname{tr} C = 0$ and that $C$ vanishes when $F$ is a multiple of the identity. Since $C$ can also be written as $(\sqrt{F^T F} - \lambda I) - \frac{1}{3}\mathbf{s}_1 I$, the remaining two invariants of $C$, $\mathbf{z}_2 = \mathcal{I}_2(C)$ and $\mathbf{z}_3 = \mathcal{I}_3(C)$, can be expressed in terms of $\mathbf{s}$, by (2.1b), and these are small for small displacements. Therefore, in terms of the variables,

$$
\begin{aligned}
(2.11) \qquad \mathbf{z}_1 &= \mathbf{r}_1 = \lambda + \tfrac{1}{3}\mathbf{s}_1, \\
\mathbf{z}_2 &= \mathcal{I}_2(C) = \mathbf{s}_2 - \tfrac{1}{3}\mathbf{s}_1^2, \\
\mathbf{z}_3 &= \mathcal{I}_3(C) = \mathbf{s}_3 - \tfrac{1}{3}\mathbf{s}_1\mathbf{s}_2 + \tfrac{2}{27}\mathbf{s}_1^3,
\end{aligned}
$$

the stored energy function has an expansion of the form

$$(2.12) \qquad \tau(\lambda, \mathbf{s}) = f(\mathbf{z}_1) + g(\mathbf{z}_1)\mathbf{z}_2 + h(\mathbf{z}_1)\mathbf{z}_3 + r(\mathbf{z}_1, \mathbf{z}_2, \mathbf{z}_3),$$

in which the remainder $r(\mathbf{z}_1, \mathbf{z}_2, \mathbf{z}_3)$ vanishes to first order in $\mathbf{z}_2$ and $\mathbf{z}_3$:

$$r(\mathbf{z}_1, 0, 0) = r_2(\mathbf{z}_1, 0, 0) = r_3(\mathbf{z}_1, 0, 0) = 0.$$

For small $\mathbf{s}$, it is clear that $f(\mathbf{z}_1)$ is the dominant term in (2.12).

The ellipticity conditions (2.5a) state that

$$(2.13a) \quad c_1^2(\lambda) = \tau_{11}(\lambda, 0) = \frac{1}{9}f''(\lambda) - \frac{2}{3}g(\lambda),$$

$$(2.13b) \quad c_2^2(\lambda) = \frac{1}{2\lambda}[\tau_1(\lambda, 0) - \lambda\tau_2(\lambda, 0)] = \frac{1}{2\lambda}[\tfrac{1}{3}f'(\lambda) - \lambda g(\lambda)].$$

Elimination of $g$ from (2.13a) and (2.13b) leads to

$$(2.14a) \qquad f''(\lambda) - \frac{2}{\lambda}f'(\lambda) = 9b(\lambda),$$

with $b(\lambda) = c_1^2(\lambda) - \frac{4}{3}c_2^2(\lambda)$, the bulk modulus which has been assumed to be positive in (2.5c). The ODE (2.14a) together with the initial condition (2.5b), which implies that $f'(1) = 0$, uniquely determines the function $f$ (up to an inessential constant) in terms of the bulk modulus

$$(2.14b) \qquad f(\lambda) = 3\int_1^\lambda \frac{\lambda^3 - y^3}{y^2} b(y)\,dy.$$

Having determined $f$ in (2.14b) and therefore also $g$ from (2.13b), the null condition (2.9c) uniquely determines the function $h$:

$$0 = \tau_{111}(\lambda, 0) = \tfrac{1}{27}f'''(\lambda) - \tfrac{4}{9}g'(\lambda) + \tfrac{4}{9}h(\lambda).$$

Starting from (2.12), the Cauchy stress for the homogeneous dilation $\varphi(x) = \lambda x$ is $\mathtt{T}(\lambda I) = (f'(\lambda)/3\lambda^2)I$, in units of force per unit area. Thus, the pressure on the sphere $|\varphi| = \lambda$ is proportional to $-f'(\lambda)/\lambda^2$. From (2.14a), positivity of the bulk modulus is seen to correspond to the monotonicity of this pressure in $\lambda$. Experimental data for rubber-like materials suggest rapid growth for the pressure as $\lambda$ becomes small [11, p. 519]. It is clear from (2.14b), that as long as $b(\lambda) \geq \mathcal{O}(\lambda^{-p})$ for $\lambda \ll 1$, with $p > 2$, the stored energy function behaves correctly near the equilibrium: $\tau(\lambda, 0) \to \infty$, as $\lambda \to 0$ and as $\lambda \to \infty$. We emphasize that the null condition places no restriction on the bulk or shear moduli of the equilibrium.

## 3. Proof of the existence theorem

3.1. *Commutation.* When we recall the definitions of the coefficients in (1.9b), (1.9c), differentiation of the relations (1.7) yields

(3.1a) $$A^{ij}_{\ell m}(\lambda) = A^{\alpha\beta}_{\delta\epsilon}(\lambda) Q_{\alpha i} Q_{\beta j} Q_{\delta \ell} Q_{\epsilon m},$$

(3.1b) $$B^{ijk}_{\ell mn}(\lambda) = B^{\alpha\beta\gamma}_{\delta\epsilon\eta}(\lambda) Q_{\alpha i} Q_{\beta j} Q_{\gamma k} Q_{\delta \ell} Q_{\epsilon m} Q_{\eta n},$$

for all proper orthogonal matrices $Q$. So $A$ and $B$ are isotropic tensors.

Consider the one-parameter family of rotations generated by the $U_\ell$ defined in (1.2b),

$$\dot{Q}_\ell(s) = U_\ell Q_\ell(s), \quad Q_\ell(0) = I.$$

If $\varphi : \mathbb{R}^3 \to \mathbb{R}^3$, define the family of simultaneous rotations

(3.2) $$T_Q \varphi(x) = Q\varphi(Q^T x).$$

Under the hypotheses (1.7), this transformation leaves the equations of motion (1.8) invariant.

From (1.9b), (3.1a) it follows that

(3.3a) $$T_{Q_\ell(s)}[Au] = A[T_{Q_\ell(s)} u],$$

and in the same way, from (1.9c), (3.1b) that

(3.3b) $$T_{Q_\ell(s)} N(u,v) = N(T_{Q_\ell(s)} u, T_{Q_\ell(s)} v),$$

In particular, we conclude from (3.3a) and (3.3b) that (3.2) also leaves the truncated system (1.9a) invariant.

The operators $\widetilde{\Omega}$ defined in (1.2a) are generated by $T_{Q_\ell(s)}$ in the sense that

(3.4a) $$\widetilde{\Omega}_\ell u = D_s T_{Q_\ell(s)} u|_{s=0}.$$

Taking the derivatives of (3.3a) and (3.3b) at $s = 0$, we see that (3.4a) gives the commutation relations

(3.4b) $$\widetilde{\Omega}_\ell A u = A \widetilde{\Omega}_\ell u,$$

and

$$\widetilde{\Omega}_\ell N(u,v) = N(\widetilde{\Omega}_\ell u, v) + N(u, \widetilde{\Omega}_\ell v). \tag{3.4c}$$

So from (3.4b) and (3.4c) any solution of the truncated linearized equation (1.9a) satisfies

$$L\widetilde{\Omega}_\ell u = N(\widetilde{\Omega}_\ell u, u) + N(u, \widetilde{\Omega}_\ell u). \tag{3.4d}$$

Next, define the one-parameter family of dilations

$$R_s u(t,x) = s^{-1} u(st, sx). \tag{3.5}$$

Since

$$R_s L u = s^{-2} L[R_s u], \tag{3.6a}$$

and

$$R_s N(u,v) = s^{-2} N(R_s u, R_s v), \tag{3.6b}$$

the scaling (3.5) leaves the truncated equations (1.9a) invariant.

The family (3.5) generates $\widetilde{S}$ defined in (1.3):

$$\widetilde{S} u(t,x) = D_s R_s u(t,x)|_{s=1}. \tag{3.7a}$$

Hence, upon differentiation of (3.6a), (3.6b), we have that

$$\widetilde{S} L u = L \widetilde{S} u - 2Lu, \tag{3.7b}$$

and

$$\widetilde{S} N(u,v) = N(\widetilde{S} u, v) + N(u, \widetilde{S} v) - 2N(u,v). \tag{3.7c}$$

By (3.7b) and (3.7c) any solution of the linearized equations (1.9a) also satisfies

$$L \widetilde{S} u = N(\widetilde{S} u, u) + N(u, \widetilde{S} u). \tag{3.7d}$$

As a consequence of (3.4d) and (3.7d), we have

PROPOSITION 3.1. *For any solution $u$ of (1.9a) in $\dot{H}^\kappa_\Gamma(T)$,*

$$L\Gamma^a u = \sum_{b+c=a} N(\Gamma^b u, \Gamma^c u), \tag{3.8}$$

*in which the sum extends over all ordered partitions of the sequence a, with $|a| \leq \kappa - 1$.*

*Note.* As a final remark in this section, we consider the commutation properties of the projections $P_\alpha$ defined in (1.5a). Since $T_Q(x/r) = x/r$, it follows that $\widetilde{\Omega}_\ell(x/r) = 0$. Thus, $\widetilde{\Omega}_\ell$ commutes with $P_\alpha$. Likewise $\partial_r(x/r) = 0$, and so $\partial_r$ also commutes with $P_\alpha$. The projections do not commute with $L$, $A$, or $\nabla$, however.

3.2. *Pointwise estimates.* The lemmas in this section use the decomposition (1.1). Despite their simplicity, the next two form the heart of the decay estimates.

LEMMA 3.1. *Let $u \in H^2_\Lambda$. For the linear operator $A$ defined in (1.9b) and its corresponding propagation speeds $c_\alpha^2$ from (2.5a),*

$$\left| P_\alpha \left[ Au(x) - c_\alpha^2 \partial_r^2 u(x) \right] \right| \leq \frac{C}{r} \left[ |\nabla \widetilde{\Omega} u(x)| + |\nabla u(x)| \right],$$

*for $\alpha = 1, 2$.*

*Proof.* Thanks to (1.1), we may write

$$\partial_\ell \partial_m u^j = \frac{x^\ell x^m}{r^2} \partial_r^2 u^j - \frac{x^\ell}{r} \partial_r \left( \frac{x}{r^2} \wedge \Omega \right)_m u^j - \left( \frac{x}{r^2} \wedge \Omega \right)_\ell \partial_m u^j,$$

from which it follows that

$$(3.9) \qquad \left| A^{ij}_{\ell m} \left( \partial_\ell \partial_m u^j - \frac{x^\ell x^m}{r^2} \partial_r^2 u^j \right) \right| \leq \frac{C}{r} \left[ |\nabla \widetilde{\Omega} u| + |\nabla u| \right],$$

by the commutation property (1.4). Recall (2.6a) which says that the symbol satisfies $A(x/r) = \sum_\alpha c_\alpha^2 P_\alpha$, so that the result is clear from (3.9). □

LEMMA 3.2. *Let $u \in \dot{H}^2_\Gamma(T)$. Then for $\alpha = 1, 2$,*

$$(3.10a) \quad |c_\alpha t - r| |P_\alpha Au(t,x)| \leq C \left[ |\nabla \Gamma u(t,x)| + |\nabla u(t,x)| + t|Lu(t,x)| \right],$$

*and*

$$(3.10b) \quad |c_\alpha t - r| |P_\alpha \partial_t \nabla u(t,x)| \leq C \left[ |\nabla \Gamma u(t,x)| + |\nabla u(t,x)| + t|Lu(t,x)| \right].$$

*Proof.* The following easily verified identities appeared in [10]:

(3.11a)
$$(c_\alpha^2 t^2 - r^2) Au(t,x) = c_\alpha^2 (t\partial_t - r\partial_r) \widetilde{S} u(t,x)$$
$$- r^2 \left[ Au(t,x) - c_\alpha^2 \partial_r^2 u(t,x) \right] - c_\alpha^2 t^2 Lu(t,x),$$

(3.11b)
$$(c_\alpha t - r) \partial_t \partial_r u(t,x) = (\partial_t - c_\alpha \partial_r) \widetilde{S} u(t,x) + \frac{(c_\alpha t - r)}{c_\alpha} Au(x,t)$$
$$+ \frac{r}{c_\alpha} \left[ Au(t,x) - c_\alpha^2 \partial_r^2 u(t,x) \right] - tLu(t,x).$$

The estimate (3.10a) is obtained from (3.11a) by dividing by $c_\alpha t + r$, applying $P_\alpha$, and using Lemma 3.1.

For the other inequality, we write

$$(c_\alpha t - r)P_\alpha \partial_t \partial_\ell u(t,x)$$
$$= (c_\alpha t - r)\left[P_\alpha \frac{x^\ell}{r}\partial_r\partial_t u(t,x) + P_\alpha \partial_t\left(\partial_\ell - \frac{x^\ell}{r}\partial_r\right)u(t,x)\right].$$

The first term satisfies the desired estimate by (3.11b) and Lemma 3.1. The second term can be estimated after first writing

$$(c_\alpha t - r)\partial_t u(x,t) = c_\alpha S u(t,x) - r[c_\alpha \partial_r + \partial_t]u(x,t).$$

Then if we apply $\partial_\ell - \frac{x^\ell}{r}\partial_r = -\left(\frac{x}{r^2}\wedge\Omega\right)_\ell$ to both sides, we find that the second term also has the right bound, since $\left(\frac{x}{r}\wedge\Omega\right)_\ell$ commutes with $c_\alpha\partial_r + \partial_t$. □

*Note.* This lemma is the only spot where the scaling operator $\widetilde{S}$ is required. A similar weighted estimate holds for $\partial_t^2 u(t,x)$, but we have no need for it.

The next result captures the manner in which the null condition will be useful in the course of the energy estimates.

PROPOSITION 3.2. *Suppose that $u,v,w \in H^2_\Lambda$. Assume that nonlinear form $N$ satisfies the null condition (2.8a). Let $\mathcal{N} = \{(\alpha,\beta,\gamma) \neq (1,1,1),(2,2,2)\}$ be the set of nonresonant indices. Then*

$$\begin{aligned}
(3.12) \quad |\langle u(x), N(v(x), w(x))\rangle| \\
\leq \frac{C}{r}|u(x)| \sum_{|a|\leq 1} \Big[|\nabla\widetilde{\Omega}^a v(x)||\nabla w(x)| + |\nabla\widetilde{\Omega}^a w(x)||\nabla v(x)| \\
+ |\nabla^2 v(x)||\widetilde{\Omega}^a w(x)| + |\nabla^2 w(x)||\widetilde{\Omega}^a v(x)|\Big] \\
+ C\sum_{\mathcal{N}} |P_\alpha u(x)|\Big[|P_\beta \nabla^2 v(x)||P_\gamma \nabla w(x)| + |P_\beta \nabla^2 w(x)||P_\gamma \nabla v(x)|\Big].
\end{aligned}$$

*Proof.* Using the projections $P_\alpha$, we write $B^{ijk}_{\ell mn} = \sum_{\alpha,\beta,\gamma} B^{\eta\mu\nu}_{\ell mn} P^{\eta i}_\alpha P^{\mu j}_\beta P^{\nu k}_\gamma$, so that

$$\begin{aligned}
(3.13) \quad \langle u, N(v,w)\rangle = \sum_{\alpha=1,2} B^{\eta\mu\nu}_{\ell mn} P^{\eta i}_\alpha P^{\mu j}_\alpha P^{\nu k}_\alpha u^i \partial_\ell(\partial_m v^j \partial_n w^k) \\
+ \sum_{\mathcal{N}} B^{\eta\mu\nu}_{\ell mn} P^{\eta i}_\alpha u^i \Big[P^{\mu j}_\beta \partial_\ell\partial_m v^j P^{\nu k}_\gamma \partial_n w^k + P^{\mu j}_\beta \partial_m v^j P^{\nu k}_\gamma \partial_\ell\partial_n w^k\Big].
\end{aligned}$$

The second group of terms in (3.13) is estimated by the second expression on the right-hand side of (3.12).

As in the proof of Lemma 3.1, the formula (1.1) enables the quantity

$$u^i \left[ \partial_\ell(\partial_m v^j \partial_n w^k) - \frac{x^\ell x^m x^n}{r^3} \partial_r(\partial_r v^j \partial_r w^k) \right]$$

to be estimated by the first group of terms on the right-hand side of (3.12). The proof of (3.12) concludes since

$$B^{\eta\mu\nu}_{\ell mn} P^{\eta i}_\alpha P^{\mu j}_\alpha P^{\nu k}_\alpha \frac{x^\ell x^m x^n}{r^3} = 0, \quad \text{for all} \quad i,j,k,\alpha$$

which follows by the null condition (2.8b) when $\alpha = 1$ and by (2.8c) for $\alpha = 2$.  □

3.3. *Sobolev inequalities.* The following Sobolev-type inequalities with weights for the most part appeared implicitly in [10].

LEMMA 3.3.  *For $u \in C_0^\infty(\mathbb{R}^3)^3$, $r = |x|$, and $\rho = |y|$,*

$$(3.14\text{a}) \qquad r^{1/2}|u(x)| \leq C \sum_{|a|\leq 1} \|\nabla \widetilde{\Omega}^a u\|_{L^2},$$

$$(3.14\text{b}) \qquad r|u(x)| \leq C \sum_{|a|\leq 1} \|\partial_r \widetilde{\Omega}^a u\|^{1/2}_{L^2(|y|\geq r)} \cdot \sum_{|a|\leq 2} \|\widetilde{\Omega}^a u\|^{1/2}_{L^2(|y|\geq r)},$$

$$(3.14\text{c}) \quad r\langle c_\alpha t - r\rangle^{1/2}|u(x)| \leq C \sum_{|a|\leq 1} \|\langle c_\alpha t - \rho\rangle \partial_r \widetilde{\Omega}^a u\|_{L^2(|y|\geq r)}$$
$$+ C \sum_{|a|\leq 2} \|\widetilde{\Omega}^a u\|_{L^2(|y|\geq r)},$$

$$(3.14\text{d}) \quad r\langle c_\alpha t - r\rangle|u(x)| \leq C \sum_{|a|\leq 1} \|\langle c_\alpha t - \rho\rangle \partial_r \widetilde{\Omega}^a u\|_{L^2(|y|\geq r)}$$
$$+ C \sum_{|a|\leq 2} \|\langle c_\alpha t - \rho\rangle \widetilde{\Omega}^a u\|_{L^2(|y|\geq r)}.$$

*Proof.* For the moment $R(r)$ will denote any smooth and positive radial function. The coordinate on the unit sphere $S^2$ will be denoted by $\omega$ and the surface measure will be $d\omega$. The proof of these inequalities begins with the bound

(3.15)
$$r^{2+\alpha} R(r)^{2+\beta} \int_{S^2} |u(r\omega)|^4 d\omega \leq C r^{2+\alpha} \int_{S^2} \int_r^\infty \left[ R(\rho)^{2+\beta} |\partial_r u(\rho\omega)||u(\rho\omega)|^3 \right.$$
$$\left. + R(\rho)^{1+\beta} |R'(\rho)||u(\rho\omega)|^4 \right] d\rho d\omega$$

$$\leq C \left( \int_{|y|\geq r} \left[ |R|^2 |\partial_r u|^2 + |R'|^2 |u|^2 \right] dy \right)^{1/2}$$

$$\times \left( \int_{|y|\geq r} |y|^{2\alpha} |R|^{2(1+\beta)} |u|^6 dy \right)^{1/2}.$$

At this point we pause to extract (3.14a) from (3.15) setting $\alpha = 0$ and $R(r) \equiv 1$. With the use of the standard estimate $\|u\|_{L^6} \leq \|\nabla u\|_{L^2}$ in $\mathbb{R}^3$, we obtain

(3.16) $$\left( r^2 \int_{S^2} |u(r\omega)|^4 d\omega \right)^{1/4} \leq C \|\nabla u\|_{L^2}.$$

Making use of another basic Sobolev inequality

(3.17) $$|u(x)| \leq C \sum_{|a|\leq 1} \|\widetilde{\Omega}^a u(r\omega)\|_{L^4(S^2)},$$

we get (3.14a) from (3.16).

We return to (3.15) in which we now take $\alpha = 2$ for the remainder of the proof. We examine the last integral in (3.15), starting with the inequality

$$\|u(r\omega)\|_{L^6(S^2)} \leq C \sum_{|a|\leq 1} \|\widetilde{\Omega}^a u(r\omega)\|_{L^2(S^2)}^{1/3} \|u(r\omega)\|_{L^4(S^2)}^{2/3}.$$

From this we get

(3.18) $$\int_{|y|\geq r} |y|^4 |R|^{2(1+\beta)} |u|^6 dy \leq C \sup_{\rho \geq r} \left( \rho^4 |R(\rho)|^{2+\beta} \int_{S^2} |u(\rho\omega)|^4 d\omega \right)$$

$$\times \sum_{|a|\leq 1} \int_{|y|\geq r} |R|^{\beta} |\widetilde{\Omega}^a u|^2 dy.$$

Putting (3.18) together with (3.15), we conclude that

(3.19)

$$\left( r^4 |R(r)|^{2+\beta} \int_{S^2} |u(r\omega)|^4 d\omega \right)^{1/4} \leq C \left( \int_{|y|\geq r} \left[ |R|^2 |\partial_r u|^2 + |R'|^2 |u|^2 \right] dy \right)^{1/4}$$

$$\times \left( \sum_{|a|\leq 1} \int_{|y|\geq r} |R|^{\beta} |\widetilde{\Omega}^a u|^2 dy \right)^{1/4}.$$

The remaining inequalities are proved by combining (3.19) with (3.17). In order to get (3.14b), take $R(r) \equiv 1$ again, while for (3.14c) and (3.14d) set $R(r) = \langle c_\alpha t - r \rangle$ (note that $|R'| \leq 1$), and $\beta = 0, 2$ respectively. □

The next result tunes these estimates to the demands of the proof by applying them to higher derivatives and also by removing the singularity near the origin.

PROPOSITION 3.3. *Let* $u \in \dot{H}^\kappa_\Gamma(T)$, *with* $\mathcal{X}_\kappa(u(t)) < \infty$.

(3.20a) $\quad \langle r \rangle^{1/2} |\Gamma^a u(t,x)| \leq C E^{1/2}_\kappa(u(t)), \qquad |a| + 2 \leq \kappa,$

(3.20b) $\quad \langle r \rangle |\partial \Gamma^a u(t,x)| \leq C E^{1/2}_\kappa(u(t)), \qquad |a| + 3 \leq \kappa,$

(3.20c) $\quad \langle r \rangle \langle c_\alpha t - r \rangle^{1/2} |P_\alpha \partial \Gamma^a u(t,x)|$
$$\leq C \left[ E^{1/2}_\kappa(u(t)) + \mathcal{X}_\kappa(u(t)) \right], \qquad |a| + 3 \leq \kappa,$$

(3.20d) $\quad \langle r \rangle \langle c_\alpha t - r \rangle |P_\alpha \partial \nabla \Gamma^a u(t,x)| \leq C \mathcal{X}_\kappa(u(t)), \qquad |a| + 4 \leq \kappa.$

*Proof.* Inequality (3.20a) and inequality (3.20b), for $r \geq 1$, follow by application of (3.14a) and (3.14b) to the indicated derivative and by the commutation property (1.4). On the other hand, inequality (3.20c) for $r \geq 1$, results from application of (3.14c) to $P_\alpha \partial \Gamma^a u$ and by the fact that $\partial_r$ and $\widetilde{\Omega}$ politely commute with $P_\alpha$, as explained at the end of Section 3.1. Similarly, (3.20d), for $r \geq 1$, arises from (3.14d). (Since (3.14c) and (3.14d) delete the origin, we need not consider the lack of smoothness of $P_\alpha u$ near the origin.) Thus, we need only prove (3.20a)–(3.20d) for $r \leq 1$.

For $r \leq 1$, (3.20b) is an immediate consequence of the Sobolev embedding

(3.21) $$H^2(\mathbb{R}^3) \subset L^\infty(\mathbb{R}^3).$$

To obtain the other inequalities, we define a cut-off function

(3.22) $$\zeta(r) = \begin{cases} 1, & \text{if } r < 1 \\ 0, & \text{if } r > 2. \end{cases}$$

Consider (3.20a) for $r \leq 1$. By (3.21),

$$|\Gamma^a u(t,x)| = \zeta(r)|\Gamma^a u(t,x)| \leq C \sum_{|b| \leq 2} \|\nabla^b(\zeta \Gamma^a u)\|_{L^2}$$
$$\leq C[E^{1/2}_\kappa(u(t)) + \|\Gamma^a u\|_{L^2(|y| \leq 2)}].$$

We finish by using (3.14a):

$$\|\Gamma^a u\|_{L^2(|y| \leq 2)} \leq C\|\rho^{1/2} \Gamma^a u\|_{L^\infty} \|\rho^{-1/2}\|_{L^2(|y| \leq 2)} \leq C E^{1/2}_\kappa(u(t)).$$

Next we turn to (3.20c), for $r \leq 1$. Note that there are constants $C < C'$ such that

(3.23a) $\quad C\zeta(r) \sum_\alpha \langle c_\alpha t - r \rangle P_\alpha \leq \zeta(r)\langle t \rangle I \leq C'\zeta(r/2) \sum_\alpha \langle c_\alpha t - r \rangle P_\alpha.$

We again make use of (3.21) and (3.23a),

$$
\begin{aligned}
(3.23\text{b}) \quad \langle c_\alpha t - r\rangle^{1/2}|P_\alpha \partial \Gamma^a u(t,x)| \\
\leq\ & C\langle t\rangle^{1/2}\zeta(r)|\partial\Gamma^a u(t,x)| \\
\leq\ & C\langle t\rangle^{1/2}\sum_{|b|\leq 2}\|\nabla^b(\zeta\partial\Gamma^a u)\|_{L^2} \\
\leq\ & C\langle t\rangle^{1/2}\sum_{|b|\leq 2}\|\partial\nabla^b\Gamma^a u\|_{L^2(|y|\leq 2)} \\
\leq\ & C\sum_\beta\sum_{|b|\leq 2}\|\langle c_\beta t-\rho\rangle^{1/2}P_\beta\partial\nabla^b\Gamma^a u\|_{L^2(|y|\leq 2)} \\
\leq\ & C\mathcal{X}_\kappa(u(t)) + C\sum_\beta \|\langle c_\beta t-\rho\rangle^{1/2}P_\beta\partial\Gamma^a u\|_{L^2(|y|\leq 2)}.
\end{aligned}
$$

To complete the proof of (3.20c), we must treat the last term above. By (3.14c),

$$
\begin{aligned}
\|\langle c_\beta t-\rho\rangle^{1/2}P_\beta\partial\Gamma^a u\|_{L^2(|y|\leq 2)} &\leq \|\rho\langle c_\beta t-\rho\rangle^{1/2}P_\beta\partial\Gamma^a u\|_{L^\infty}\|\rho^{-1}\|_{L^2(|y|\leq 2)} \\
&\leq C[E_\kappa^{1/2}(u(t)) + \mathcal{X}_\kappa(u(t))].
\end{aligned}
$$

The proof of (3.20d) for $r\leq 1$ follows the same lines as (3.23b), except that the last step above is no longer necessary. □

*Note.* In the Lorentz invariant case, stronger versions of these inequalities hold which avoid the weighted norm $\mathcal{X}_\kappa$ altogether.

3.4. *Weighted $L^2$-estimates.* The goal of this section is to show that the weighted norm $\mathcal{X}_\kappa$ is controlled by the energy $E_\kappa^{1/2}$, for small solutions of (1.9a).

LEMMA 3.4. *Let $u \in \dot{H}_\Gamma^2(T)$. Then*

$$(3.24) \quad \mathcal{X}_2(u(t)) \leq C\left[E_2^{1/2}(u(t)) + t\|Lu(t)\|_{L^2}\right].$$

*Proof.* Recall definition (1.6):

$$(3.25) \quad \mathcal{X}_2(u(t)) = \sum_\alpha \|\langle c_\alpha t - r\rangle P_\alpha \partial \nabla u(t)\|_{L^2}.$$

In view of (3.10b), it is enough to consider the terms with two spatial derivatives.

We introduce a smooth cut-off function along $r = \langle c_2 t\rangle/2$, by defining

$$(3.26) \quad \psi(t,r) = \zeta(4r/\langle c_2 t\rangle),$$

using $\zeta$ as defined in (3.22). On the support of $\psi$, the ratio $\langle c_\alpha t - r\rangle/\langle t\rangle$ is uniformly bounded above and below, and so we have

$$(3.27) \quad C\psi(t,r)\langle t\rangle I \leq \psi(t,r)\sum_\alpha \langle c_\alpha t - r\rangle P_\alpha \leq C'\psi(t,r)\langle t\rangle I.$$

Thus by (3.27), the terms of interest have the estimate

$$\sum_\alpha \|\langle c_\alpha t - r\rangle P_\alpha \nabla^2 u(t)\|_{L^2} \tag{3.28}$$

$$\leq \sum_\alpha \|\langle c_\alpha t - r\rangle \psi P_\alpha \nabla^2 u(t)\|_{L^2}$$

$$+ \sum_\alpha \|\langle c_\alpha t - r\rangle (1-\psi) P_\alpha \nabla^2 u(t)\|_{L^2}$$

$$\leq \langle t\rangle \|\psi \nabla^2 u(t)\|_{L^2} + \sum_\alpha \|\langle c_\alpha t - r\rangle(1-\psi) P_\alpha \nabla^2 u(t)\|_{L^2}.$$

First we will show that the integral over the support of $\psi$ is controlled by the elliptic operator $A$. From the explicit formula (2.6b), we have

$$\|\psi A u(t)\|_{L^2}^2 = c_2^4 \int \psi^2 |\Delta u(t)|^2 dx \tag{3.29a}$$

$$+ 2c_2^2(c_1^2 - c_2^2)\int \psi^2 \langle \Delta u(t), \nabla(\nabla \cdot u(t))\rangle dx$$

$$+ (c_1^2 - c_2^2)^2 \int \psi^2 |\nabla(\nabla \cdot u(t))|^2 dx.$$

Since the first order spatial derivatives of $\psi$ are $\mathcal{O}(\langle t\rangle^{-1})$, integration by parts shows that

$$\int \psi^2 |\Delta u(t)|^2 dx \tag{3.29b}$$

$$\geq \sum_{\ell,m} \|\psi \partial_\ell \partial_m u(t)\|_{L^2}^2 - C\langle t\rangle^{-1}\|\psi \nabla^2 u(t)\|_{L^2}\|\nabla u(t)\|_{L^2}.$$

Likewise,

$$\int \psi^2 \langle \Delta u(t), \nabla(\nabla \cdot u(t))\rangle dx \tag{3.29c}$$

$$\geq \|\psi \nabla(\nabla \cdot u(t))\|_{L^2}^2 - C\langle t\rangle^{-1}\|\psi \nabla^2 u(t)\|_{L^2}\|\nabla u(t)\|_{L^2}.$$

From (3.29a), (3.29b), and (3.29c), we obtain, using (3.27) again,

(3.29d)

$$\sum_{\ell,m} c_2^2 \langle t\rangle \|\psi \partial_\ell \partial_m u(t)\|_{L^2} \leq C\Big[E_2^{1/2}(u(t)) + \langle t\rangle \|\psi A u(t)\|_{L^2}\Big]$$

$$\leq C\Big[E_2^{1/2}(u(t)) + \sum_\alpha \|\langle c_\alpha t - r\rangle \psi P_\alpha A u(t)\|_{L^2}\Big].$$

On the support of $1-\psi$, we note that the quantity $\langle c_\alpha t - r\rangle/r$ is bounded. So for the integrals over the support of $1-\psi$, we use first (1.1) to pass from an arbitrary second spatial derivative to radial derivatives, and then Lemma 3.1 to pass to $A$,

$$
\begin{aligned}
(3.30) \quad \|\langle c_\alpha t - r\rangle(1-\psi)P_\alpha \partial_\ell \partial_m u(t)\|_{L^2} \\
\leq C\Big[E_2^{1/2}(u(t)) + \|\langle c_\alpha t - r\rangle(1-\psi)P_\alpha \partial_r^2 u(t)\|_{L^2}\Big] \\
\leq C\Big[E_2^{1/2}(u(t)) + c_\alpha^{-2}\|\langle c_\alpha t - r\rangle(1-\psi)P_\alpha A u(t)\|_{L^2}\Big].
\end{aligned}
$$

As a consequence of (3.25), (3.28), (3.29d), and (3.30), we have that

$$
\begin{aligned}
\mathcal{X}_2(u(t)) \leq C\Big[&E_2^{1/2}(u(t)) + \sum_\alpha \|\langle c_\alpha t - r\rangle P_\alpha A u(t)\|_{L^2} \\
&+ \sum_\alpha \|\langle c_\alpha t - r\rangle P_\alpha \partial_t \nabla u(t)\|_{L^2}\Big],
\end{aligned}
$$

and so the lemma follows now from Lemma 3.2. □

Next we estimate the nonlinear terms on the right-hand side of (3.24).

LEMMA 3.5. *Let $u \in \dot{H}_\Gamma^\mu(T)$ be a solution of (1.9a). Then*

$$\mathcal{X}_\mu(u(t)) \leq C\Big[E_\mu^{1/2}(u(t)) + \mathcal{X}_{\mu'}(u(t))E_\mu^{1/2}(u(t)) + \mathcal{X}_\mu(u(t))E_{\mu'}^{1/2}(u(t))\Big],$$

*with $\mu' = \left[\frac{\mu-1}{2}\right] + 3$.*

*Proof.* First applying Lemma 3.4 to $\Gamma^a u$ and summing over $|a| \leq \mu - 2$, we have that

$$(3.31) \qquad \mathcal{X}_\mu(u(t)) \leq CE_\mu^{1/2}(u(t)) + C\sum_{|a|\leq\mu-2} t\|L\Gamma^a u(t)\|_{L^2}.$$

Because of Proposition 3.1, we need only to estimate terms of the form

$$t\|\nabla^2 \Gamma^b u \nabla \Gamma^c u\|_{L^2},$$

with $|b+c| = |a| \leq \mu-2$. Plainly either $|b| \leq m-1$ or $|c| \leq m$ where $m = [\frac{\mu-1}{2}]$.

We note that the quantities $\langle t\rangle^{-1}\langle r\rangle\langle c_\alpha t - r\rangle$ are uniformly bounded below, so that

$$(3.32) \qquad I \leq C\langle t\rangle^{-1}\langle r\rangle \sum_\alpha \langle c_\alpha t - r\rangle P_\alpha.$$

Thus we have the estimate

$$t\|\nabla^2 \Gamma^b u \nabla \Gamma^c u\|_{L^2} \leq C\sum_\alpha \|\langle r\rangle\langle c_\alpha t - r\rangle P_\alpha \nabla^2 \Gamma^b u \nabla \Gamma^c u\|_{L^2}$$

$$\leq C \begin{cases} \sum_\alpha \|\langle r\rangle\langle c_\alpha t - r\rangle P_\alpha \nabla^2 \Gamma^b u\|_{L^\infty}\|\nabla \Gamma^c u\|_{L^2}, & |b| \leq m-1 \\ \sum_\alpha \|\langle c_\alpha t - r\rangle P_\alpha \nabla^2 \Gamma^b u\|_{L^2}\|\langle r\rangle \nabla \Gamma^c u\|_{L^\infty}, & |c| \leq m. \end{cases}$$

In the first case, using (3.20d) we get the upper bound

$$\mathcal{X}_{m+3}(u(t))E_{\mu-1}^{1/2}(u(t)),$$

and in the second case, by (3.20b) we get the upper bound

$$\mathcal{X}_\mu(u(t))E_{m+3}^{1/2}(u(t)). \qquad \square$$

The following bootstrap argument completes the estimation of $\mathcal{X}_\kappa$. Notice that the smallness condition is imposed only for $\mathcal{X}_{\kappa-2}$.

PROPOSITION 3.4. *Let $u \in \dot{H}_\Gamma^\kappa(T)$, $\kappa \geq 8$, be a solution of (1.9a). If $E_\mu(u(t))$, $\mu = \kappa - 2$, remains sufficiently small for $0 \leq t < T$, then*

(3.33a) $\qquad \mathcal{X}_\mu(u(t)) \leq CE_\mu^{1/2}(u(t)),$

(3.33b) $\qquad \mathcal{X}_\kappa(u(t)) \leq CE_\kappa^{1/2}(u(t))[1 + E_\mu^{1/2}(u(t))].$

*Proof.* Since we have $\mu \geq 6$, it follows that $\mu' = [\frac{\mu-1}{2}] + 3 \leq \mu$. Thus, by Lemma 3.5,

$$\mathcal{X}_\mu(u(t)) \leq C\left[E_\mu^{1/2}(u(t)) + E_\mu^{1/2}(u(t))\mathcal{X}_\mu(u(t))\right],$$

and so we see that for $E_\mu^{1/2}(u(t))$ small enough, the bound (3.33a) holds.

Since $\kappa \geq 8$, it is true that $\kappa' = [\frac{\kappa-1}{2}] + 3 \leq \kappa - 2 = \mu$. So again by Lemma 3.5, we may write

$$\mathcal{X}_\kappa(u(t)) \leq C\left[E_\kappa^{1/2}(u(t)) + \mathcal{X}_\mu(u(t))E_\kappa^{1/2}(u(t)) + \mathcal{X}_\kappa(u(t))E_\mu^{1/2}(u(t))\right].$$

If $E_\mu^{1/2}(u(t))$ is small, then this implies that

$$\mathcal{X}_\kappa(u(t)) \leq CE_\kappa^{1/2}(u(t))[1 + \mathcal{X}_\mu(u(t))].$$

Thus, we obtain (3.33b) from this and (3.33a). $\qquad \square$

3.5. *Energy estimates.* Suppose that $u(t) \in \dot{H}_\Gamma^\kappa(T)$ is a local solution of (1.9a). Let $\mu = \kappa - 2$, and assume that $\varepsilon(\lambda)$, and hence $E_\mu(u(0))$, is small enough in (1.12) so that Proposition 3.4 holds at least for a short time. We suppose that $T$ is the largest such time with this property.

Proposition (3.1) and the energy method yield

(3.34) $\qquad E_\kappa'(u(t)) = \sum_{|a|\leq \kappa-1} \sum_{b+c=a} \int \langle \partial_t \Gamma^a u, N(\Gamma^b u, \Gamma^c u)\rangle dx.$

Terms in (3.34) with $b = a$ or $c = a$, $|a| = \kappa - 1$, require special attention. Thanks to (1.10) they can be written as follows:

$$
\begin{aligned}
(3.35) \quad \int \langle \partial_t \Gamma^a u, N(\Gamma^a u, u) \rangle dx &= B^{ijk}_{\ell mn} \int \partial_t \Gamma^a u^i \partial_\ell (\partial_m \Gamma^a u^j \partial_n u^k) dx \\
&= -B^{ijk}_{\ell mn} \int \partial_t \partial_\ell \Gamma^a u^i \partial_m \Gamma^a u^j \partial_n u^k dx \\
&= -\frac{1}{2} \frac{d}{dt} B^{ijk}_{\ell mn} \int \partial_\ell \Gamma^a u^i \partial_m \Gamma^a u^j \partial_n u^k dx \\
&\quad + \frac{1}{2} B^{ijk}_{\ell mn} \int \partial_\ell \Gamma^a u^i \partial_m \Gamma^a u^j \partial_t \partial_n u^k dx.
\end{aligned}
$$

If we set

$$\widetilde{N}(u, v, w) = B^{ijk}_{\ell mn} \partial_\ell u^i \partial_m v^j \partial_n w^k,$$

and

$$(3.36a) \quad \widetilde{E}_\kappa(u(t)) = E_\kappa(u(t)) + \sum_{|a| \leq \kappa - 1} \int \widetilde{N}(\Gamma^a u, \Gamma^a u, u) dx,$$

then (1.11), (3.34) and (3.35) show that

$$(3.36b) \quad \widetilde{E}'_\kappa(u(t)) = \sum_{|a| \leq \kappa - 1} \sum_{\substack{b+c=a \\ b,c \neq a}} \int \langle \partial_t \Gamma^a u, N(\Gamma^b u, \Gamma^c u) \rangle dx$$

$$+ \sum_{|a| \leq \kappa - 1} \int \widetilde{N}(\partial_t u, \Gamma^a u, \Gamma^a u) dx.$$

From (3.36b), we get

$$\widetilde{E}'_\kappa(u(t)) \leq C \sum_{|a| \leq \kappa - 1} \sum_{\substack{b+c=a \\ b \neq a}} \|\partial \Gamma^a u\|_{L^2} \|\partial \nabla \Gamma^b u \nabla \Gamma^c u\|_{L^2}.$$

In the sum, we may assume that either $|b| \leq d - 1$ or $|c| \leq d$, where $d = [\kappa/2]$. As in the proof of Lemma 3.5, we get

$$\widetilde{E}'_\kappa(u(t)) \leq C\langle t \rangle^{-1} E^{1/2}_\kappa(u(t)) \begin{cases} \mathcal{X}_\kappa(u(t)) E^{1/2}_{d+3}(u(t)), & |c| \leq d \\ \mathcal{X}_{d+3}(u(t)) E^{1/2}_\kappa(u(t)), & |b| \leq d - 1. \end{cases}$$

Since $\kappa \geq 9$, we have that $d + 3 \leq \mu = \kappa - 2$, and so by Proposition 3.4 we obtain

$$(3.37) \quad \widetilde{E}'_\kappa(u(t)) \leq C\langle t \rangle^{-1} E^{1/2}_\mu(u(t)) E_\kappa(u(t)).$$

(This estimate implies almost global existence, [5], [10].)

The next and decisive step is to show that the lower energy $E_\mu(u(t))$ remains small. Here finally the null condition enters the fray. From (3.36a), (3.36b), we have

$$(3.38a) \quad \widetilde{E}'_\mu(u(t)) = \sum_{|a|\leq\mu-1} \sum_{\substack{b+c=a \\ b,c\neq a}} \left[ \int \psi \langle \partial_t \Gamma^a u, N(\Gamma^b u, \Gamma^c u)\rangle dx \right.$$
$$\left. + \int (1-\psi)\langle \partial_t \Gamma^a u, N(\Gamma^b u, \Gamma^c u)\rangle dx \right]$$
$$+ \sum_{|a|\leq\mu-1} \left[ \int \psi \widetilde{N}(\partial_t u, \Gamma^a u, \Gamma^a u) dx \right.$$
$$\left. + \int (1-\psi)\widetilde{N}(\partial_t u, \Gamma^a u, \Gamma^a u) dx \right]$$

in which we again use the cut-off defined in (3.26). We shall focus on the terms involving the quadratic form $N$ since the ones with the trilinear form $\widetilde{N}$ are similar.

The integrals in (3.38a) over the support of $\psi$ where, it will be recalled the ratio $\langle t\rangle^{-1}\langle c_\alpha t - r\rangle$ is bounded below, are estimated using (3.20c) and Proposition 3.4

$$(3.38b) \quad \int \psi \langle \partial_t \Gamma^a u, N(\Gamma^b u, \Gamma^c u)\rangle dx$$
$$\leq C\|\partial_t \Gamma^a u\|_{L^2}\|\psi \nabla^2 \Gamma^b u \nabla \Gamma^c u\|_{L^2}$$
$$\leq C\langle t\rangle^{-3/2} E_\mu^{1/2}(u(t))$$
$$\times \sum_{\alpha,\beta} \|\langle c_\alpha t - r\rangle P_\alpha \nabla^2 \Gamma^b u\|_{L^2} \|\langle c_\beta t - r\rangle^{1/2} P_\beta \nabla \Gamma^c u\|_{L^\infty}$$
$$\leq C\langle t\rangle^{-3/2} E_\mu^{1/2}(u(t)) \mathcal{X}_{|b|+2}(u(t)) \left[ E_{|c|+3}^{1/2}(u(t)) + \mathcal{X}_{|c|+3}(u(t)) \right]$$
$$\leq C\langle t\rangle^{-3/2} E_\mu^{1/2}(u(t)) E_{|b|+2}^{1/2}(u(t)) E_{|c|+3}^{1/2}(u(t))$$
$$\leq C\langle t\rangle^{-3/2} E_\mu(u(t)) E_\kappa^{1/2}(u(t)),$$

since $|b|+2 \leq \mu$ and $|c|+3 \leq \mu+1 < \kappa$. Observe that the lack of a weighted estimate for the derivative $\nabla \Gamma^c u$ in $L^2$ forces us to estimate it in $L^\infty$, here and later on as well.

To the second group of integrals in (3.38a), we apply Proposition 3.2. This leads to a bound of the form

(3.38c) $\quad C \sum_{|a|\leq\mu-1} \sum_{\substack{b+c=a \\ b,c\neq a}} \Big\{ \|\partial_t \Gamma^a u\|_{L^2} \Big[ \|r^{-1}(1-\psi)\nabla\Gamma^{b+1}u\,\nabla\Gamma^c u\|_{L^2}$

$$+ \|r^{-1}(1-\psi)\nabla^2\Gamma^b u\,\Gamma^{c+1}u\|_{L^2} \Big]$$

$$+ \sum_{\mathcal{N}} \int |P_\alpha \partial_t \Gamma^a u||P_\beta \nabla^2 \Gamma^b u||P_\gamma \nabla \Gamma^c u| dx \Big\}.$$

Since $r \geq \langle t \rangle$ on the support of $1 - \psi$, the first group of terms in (3.38c) is bounded by

(3.38d) $\quad C\langle t\rangle^{-3/2} \sum_{|a|\leq\mu-1} \sum_{\substack{b+c=a \\ b,c\neq a}} \|\partial_t \Gamma^a u\|_{L^2}$

$$\times \Big[ \|\nabla\Gamma^{b+1}u\|_{L^2}\|r\nabla\Gamma^c u\|_{L^\infty} + \|\nabla^2\Gamma^b u\|_{L^2}\|r^{1/2}\Gamma^{c+1}u\|_{L^\infty} \Big].$$

Since $|b| + 2 \leq \mu$ and $|c| + 3 \leq \mu + 1 < \kappa$, we see by (3.20a) and (3.20b) that this expression is also bounded by

$$C\langle t\rangle^{-3/2} E_\mu(u(t)) E_\kappa^{1/2}(u(t)).$$

Now we consider the nonresonant terms in (3.38c), which will be estimated using (3.20c) and Proposition 3.4. Note that

$$1 \leq C\langle t\rangle^{-3/2}\langle r\rangle\langle c_\beta t - r\rangle\langle c_\gamma t - r\rangle^{1/2}, \quad \text{if} \quad \beta \neq \gamma.$$

Thus, if $\beta \neq \gamma$, we have

(3.38e)

$$\int |P_\alpha \partial_t \Gamma^a u||P_\beta \nabla^2 \Gamma^b u||P_\gamma \nabla \Gamma^c u| dx$$

$$\leq C\langle t\rangle^{-\frac{3}{2}} \int |P_\alpha \partial_t \Gamma^a u|\langle c_\beta t - r\rangle|P_\beta \nabla^2 \Gamma^b u|\langle r\rangle\langle c_\gamma t - r\rangle^{\frac{1}{2}}|P_\gamma \nabla\Gamma^c u| dx$$

$$\leq C\langle t\rangle^{-\frac{3}{2}} \|\partial_t \Gamma^a u\|_{L^2} \|\langle c_\beta t - r\rangle P_\beta \nabla^2\Gamma^b u\|_{L^2} \|\langle r\rangle\langle c_\gamma t - r\rangle^{\frac{1}{2}} P_\gamma \nabla\Gamma^c u\|_{L^\infty}$$

$$\leq C\langle t\rangle^{-3/2} E_\mu(u(t)) E_\kappa^{1/2}(u(t)),$$

since $|a| + 1 \leq \mu$, $|b| + 2 \leq \mu$, and $|c| + 3 < \kappa$.

In the remaining cases where $\alpha \neq \beta$, the weight $\langle r\rangle\langle c_\alpha t - r\rangle^{1/2}$ can be used in front of $\partial_t \Gamma^a u$ and estimated in $L^\infty$ by $E_\kappa^{1/2}(u(t))$, with the help of (3.20c).

Putting everything together (3.38a)–(3.38e), we have shown that

(3.39) $\quad \widetilde{E}'_\mu(u(t)) \leq C\langle t\rangle^{-3/2} E_\mu(u(t)) E_\kappa^{1/2}(u(t)).$

The smallness of the energy implies, in turn, the smallness of $|\nabla u|$ by (3.20b). Thus, by possible further restriction of $\varepsilon(\lambda)$, the modified energy in (3.36a) is equivalent to the standard one, for small solutions, and from (3.37), (3.39) we arrive at the coupled pair of differential inequalities

$$\begin{aligned}\widetilde{E}'_\kappa(u(t)) &\leq C\langle t\rangle^{-1}\widetilde{E}_\mu^{1/2}(u(t))\widetilde{E}_\kappa(u(t)),\\ \widetilde{E}'_\mu(u(t)) &\leq C\langle t\rangle^{-3/2}\widetilde{E}_\mu(u(t))\widetilde{E}_\kappa^{1/2}(u(t)).\end{aligned}$$

Now, the global bounds given in Theorem 1.1 are immediate, provided $\varepsilon(\lambda)$ is small enough.

*Acknowledgment.* The author gratefully acknowledges many informative discussions with A. Shadi Tahvildar-Zadeh leading to the preparation of this manuscript.

UNIVERSITY OF CALIFORNIA, SANTA BARBARA, CA
*E-mail address*: sideris@math.ucsb.edu